\documentclass[12pt]{article}

\usepackage[LGR,T1]{fontenc}
\usepackage[greek,english]{babel}
\usepackage{graphicx,amssymb,amsmath,amsthm}
\usepackage{hyperref,url}
\hypersetup{
  colorlinks=true,
  linkcolor=blue,
  citecolor=blue,
  urlcolor=blue
}
\usepackage{appendix}
\usepackage{booktabs}
\usepackage{verbatim}
\usepackage{geometry}

\newcommand{\nocopyright}{
No Copyright\thanks{
\href{https://creativecommons.org/publicdomain/zero/1.0/}{CC0}:
to the extent possible under law, Peter Doyle has waived all copyright
and related and neighboring rights to this work.
}}
\title{Neoplatonic solids}
\author{Peter Doyle \and Matthew Ellison}
\date{Version dated 2026-07-27
\\ \nocopyright
}

\newtheorem{theorem}{Theorem}
\newtheorem{proposition}[theorem]{Proposition}

\newtheorem*{hypothesis}{Hypothesis}

\newtheorem{corollary}[theorem]{Corollary}
\newtheorem*{conjecture}{Conjecture}

\newcommand{\proofstart}{{\noindent \bf Proof.\ }}

\newcommand{\proofend}{$\quad \qed$}
\newcommand{\goesto}{\rightarrow}

\newcommand{\beginfigurepos}{\begin{figure}[htbp]}
\newcommand{\fig}[3]{
\beginfigurepos
\includegraphics[width=370pt]{figures/#1.pdf}
\caption{#3}
\label{#2}
\end{figure}
}
\newcommand{\figsize}[4]{
\beginfigurepos
\centerline{
\includegraphics[width=#1]{figures/#2.pdf}
}
\caption{#4}
\label{#3}
\end{figure}
}

\newcommand{\bigfig}[3]{
\begin{figure}[h!]
\centerline{
\includegraphics[width=\textwidth]{figures/#1.pdf}
}
\caption{#3}
\label{#2}
\end{figure}
}

\newcommand{\sect}[1]{\S\ #1}

\makeatletter
\newcommand{\xRrightarrow}[2][]{\ext@arrow 0359\Rrightarrowfill@{#1}{#2}}
\newcommand{\Rrightarrowfill@}{\arrowfill@\equiv\equiv\Rrightarrow}
\newcommand{\xLleftarrow}[2][]{\ext@arrow 3095\Lleftarrowfill@{#1}{#2}}
\newcommand{\Lleftarrowfill@}{\arrowfill@\Lleftarrow\equiv\equiv}
\newcommand{\xLleftRrightarrow}[2][]{\ext@arrow 3399\LleftRrightarrowfill@{#1}{#2}}
\newcommand{\LleftRrightarrowfill@}{\arrowfill@\Lleftarrow\equiv\Rrightarrow}
\makeatother

\newcommand{\kw}{\mathbf}

\newcommand{\len}{\kw{len}}

\newcommand{\puff}{\Delta}
\newcommand{\hpuff}{\Delta}

\begin{document}

\maketitle

\begin{abstract}
A \emph{6-net} is a simplicial triangulation of the $2$-sphere with maximum degree $\leq 6$.
Experiments suggest that every $6$-net admits a unique realization
as an undented Euclidean polyhedron built from unit equilateral triangles,
and a unique realization as an ideal equilateral hyperbolic polyhedron.
We call these \emph{neoplatonic solids} and \emph{ideal neoplatonics}.

A net is \emph{prime} if every 3-cycle bounds a face.
A computer-assisted proof shows
that every prime $6$-net with $v \leq 50$
has a unique realization as a convex ideal neoplatonic.
Numerical homotopy from this realization yields
an approximate Euclidean neoplatonic,
and a computer-assisted proof shows that a true Euclidean neoplatonic 
lies nearby, though we do not prove uniqueness.
Using the separating-triangle decomposition,
we extend Euclidean existence to all $10{,}412{,}340$ $6$-nets
with $v\leq50$, counted up to combinatorial isomorphism.
\end{abstract}
\begin{center}
{\setlength{\fboxsep}{6pt}\fbox{\textgreek{Jeait'htw| >Ajhna'iw|}}}
\end{center}

\section{Atlas}
This paper should be read alongside an
\href{https://math.dartmouth.edu/~doyle/docs/atlas/}{online atlas of neoplatonic solids}
\cite{atlas}.

\section{Introduction}

A \emph{6-net} is a simplicial triangulation of the $2$-sphere with maximum degree
at most $6$.
(Cf.\ Thurston \cite{thurston:shapes}.)
A net is prime if every 3-cycle bounds a face.
A 6-net is prime just if it is the tetrahedron,
or if it has no vertex of degree 3
and is not a stack of two or more octahedra.

Experiments suggest:
\begin{conjecture}[Neoplatonic]
Any $6$-net $\sigma$ has a realization,
unique up to isometry,
as an undented Euclidean polyhedron $\Delta(\sigma)$
built from equilateral triangles of side length $1$.
\end{conjecture}

A polyhedron is \emph{undented} if every vertex has a local support
plane:
From outside the object you can touch any vertex with your palm
if your hand is small enough.
Equivalently,
at every vertex the sum of the exterior dihedral angles is $\geq 0$.
By Gauss--Bonnet, this guarantees that when you look from a vertex into the interior of the solid,
its immediate neighborhood occupies at most half of the visual sphere.
The atlas gives examples of dented polyhedra,
and of polyhedra with vertices that have a local but no global support plane.

We call these realizations of $6$-nets
\emph{neoplatonic solids},
as they follow in the footsteps of the tetrahedron, octahedron,
and icosahedron.
We call them `solids' out of respect for tradition,
though we regard a polyhedron as a $2$-dimensional object
rather than the region of space that it bounds.
Note that the cube and dodecahedron don't count as neoplatonics:
We may view them as subsidiary to their duals.

The convex gallery in the atlas
starts with the eight neoplatonics that are convex without coplanar faces.
They were identified by Rausenberger
\cite{rausenberger}, and later by Freudenthal and van der Waerden
\cite{freudenthal:euclid}.
(One of these eight, the triangular bipyramid, is non-prime.)
Allowing coplanar faces gives infinitely many more convex neoplatonics.
Note that we allow as a limiting case `pancakes':
polyhedra with empty interior obtained by doubling
a portion of the triangular lattice.
Pancakes are limits of embedded realizations of the same triangulated sphere.
Here the simplicial hypothesis matters:
a vertex of degree $2$ would permit bending along the opposite $2$-cycle
without creating dents, contradicting uniqueness.

Beyond these convex examples lies a world of neoplatonic solids
that are not convex but still undented.
We like to think that Theaetetus,
the discoverer of the icosahedron,
and Euclid, who constructs the icosahedron in
Book XIII, Proposition 16, of the \emph{Elements}
\cite{euclid:peyrard3},
would be pleased to see
the tetrahedron, octahedron, and icosahedron
as the eldest siblings
in a distinguished infinite family.

\begin{proposition} \label{euclid}
There are 8,239,684 prime 6-nets with $v \le 50$.
Every one has
a neoplatonic realization.
\end{proposition}

\proofstart
We generate prime 6-nets using a variant of the method of Goedgebeur
\cite{ghent}.
To prove they have neoplatonic realizations we use the method of Ellison
\cite{ellison:embed}.
We compute an approximate neoplatonic,
then use an effective version of the inverse function theorem
combined with arithmetic with guaranteed error bounds to prove that
an exact neoplatonic lies nearby.
In cases where the limit lies flat at one or more degree-6 vertices,
some of the equations specifying edge lengths get replaced by coplanarity conditions.

The software used for this computation is archived in
\cite{euclidean:software}.
\proofend

\section{Non-prime neoplatonics}

The atlas includes a gallery of non-prime neoplatonics.
Cutting along all separating triangles decomposes a non-prime net
into prime factors.
For $v\leq50$, the non-prime nets fall into three disjoint classes.
\begin{enumerate}
\item
All prime factors are tetrahedra.
With $n$ tetrahedral factors, the numbers for $n=2,3,4,5$ are
respectively $1,1,3,3$.
For each $6\leq n\leq47$ there is one net,
the helical stack called `3NA' in the atlas.
\item
There are two or more non-tetrahedral factors.
They are octahedra arranged in a stack,
with a tetrahedral cap of depth $0$, $1$, or $2$ at each end.
\item
There is exactly one non-tetrahedral prime factor,
with one to four tetrahedra attached.
\end{enumerate}

The first class contains
\[
1+1+3+3+42=50
\]
nets with $v\leq50$.
For the second class, let $m\geq2$ be the number of octahedral factors,
and let $t$ be the total depth of the two caps.
Then
\[
v=3+3m+t.
\]
For $t=0,1,2,3,4$, the numbers of cap types are respectively
$1,1,2,1,2$.
The final $2$ counts the two inequivalent relative attachments
when both caps have depth $2$.
Thus the second class contains
\[
14+14+2\mathbin{\cdot}14+13+2\mathbin{\cdot}13=95
\]
nets with $v\leq50$.

The first two classes have direct realizations by gluing regular cells.
To check them without roundoff, scale so that the squared edge length is $8$.
If $g$ is the centroid of an outward-oriented exposed face $(a,b,c)$ and
$n=(b-a)\mathbin{\times}(c-a)$,
the new tetrahedral vertex is $g+n/3$,
while the center of an attached octahedron is $g+n/6$;
its three new vertices are the reflections of $a,b,c$ through that center.
Thus all coordinates are rational.
An exact check of all $145$ nets proves that boundary triangles meet
only in their common simplices and that every vertex has a local support plane.

For the third class,
extend the certified intervals for the exterior dihedral angles
across the attached tetrahedra.
After identifying combinatorial duplicates,
this yields $2{,}172{,}511$ distinct nets with $v\leq50$.
An interval embeddedness check proves that each has
an exact embedded realization.

These attachments preserve undentedness.
Indeed, let $A\leq2\pi$ be the area inside the spherical link
of a core vertex of degree $d$.
The link has perimeter $d\pi/3$, so spherical isoperimetry gives
\[
2\pi-A\ \geq\ \frac{\pi}{3}\sqrt{36-d^2}.
\]
Each incident tetrahedral attachment adds
$E=3\arccos(1/3)-\pi<\pi/2$ to $A$,
and the degree bound allows at most $6-d$ such attachments.
For $d=4,5$ the displayed lower bound exceeds respectively
$2E,E$; for $d=6$ no attachment is possible.
New vertices have link area $E$ or $2E$.
Thus every attached realization is undented.
The software and check outputs for these computations are archived in
\cite{euclidean:software}.

\begin{corollary}\label{all-six-nets}
Every $6$-net with $v\leq50$ has a neoplatonic realization.
Up to combinatorial isomorphism, there are $10{,}412{,}340$ such nets:
$8{,}239{,}684$ prime and $2{,}172{,}656$ non-prime.
\end{corollary}

\proofstart
The proposition handles the prime case.
The three classes above contain
\[
50+95+2{,}172{,}511=2{,}172{,}656
\]
non-prime nets, and the preceding constructions realize each of them.
For every $v\leq50$, these counts together with the prime count
sum to the census of $\{3,4,5,6\}$-triangulations in
\cite[Tables~1 and~2]{brinkmann:restricted}.
\proofend

\section{Ideal neoplatonics}

Experiments suggest
that Euclidean neoplatonics are the shriveled remains of
equilateral ideal hyperbolic polyhedra,
which we call \emph{ideal neoplatonics}.
By \emph{equilateral} we mean that the ideal triangle faces
meet midpoint-to-midpoint,
where the midpoint of an infinite edge is the foot of the perpendicular dropped
from the opposite vertex.
(Equivalently, after adjacent faces are unfolded into a common plane,
their incircles are tangent.)

\begin{conjecture}[Ideal neoplatonic]
Every $6$-net $\sigma$ has
a realization $\hpuff(\sigma,\infty)$,
unique up to isometry,
as an ideal neoplatonic.
\end{conjecture}

\begin{conjecture}[Ideal prime]
If the $6$-net $\sigma$ is prime,
its ideal neoplatonic realization $\hpuff(\sigma,\infty)$
is convex.
\end{conjecture}

\begin{proposition}
Every prime $6$-net with $v \leq 50$
has a convex ideal neoplatonic realization,
unique up to isometry.
\end{proposition}

\proofstart
We use the approach of Bobenko, Pinkall, and Springborn
\cite{bobenko:discrete}.
Existence is proved using a Perron-type method of super- and subsolutions
for the equations of Springborn
\cite{springborn:ideal};
for details see \sect{\ref{sec:real}}.
Uniqueness comes from Rivin
\cite{rivin:one},
who proves that any net can have at most one
convex equilateral ideal realization.
\proofend
\par

Rivin
\cite{rivin:one}
shows that a convex ideal equilateral polyhedron maximizes volume
among convex ideal realizations of its net.
It is natural to ask whether convexity of the competitors can be dropped.
For an oriented ideal cycle, we use algebraic volume.

\begin{conjecture}[Volume max]
For a prime $6$-net $\sigma$,
its ideal neoplatonic realization maximizes volume
among all ideal geodesic $2$-cycles with combinatorics $\sigma$.
\end{conjecture}

The corresponding assertion for non-prime $6$-nets follows by splitting
along separating triangles and maximizing the resulting cycles independently.
Without the degree bound, maximizing volume is more problematic.
A regular ideal $12$-gonal bipyramid maximizes volume among convex
realizations, but twelve regular ideal tetrahedra can wind twice around
a common axis, producing a branched realization of maximal possible volume.

\section{Neoplatonic homotopy}

\begin{conjecture}[Neoplatonic homotopy]
Every $6$-net $\sigma$ has a family of realizations $\hpuff(\sigma,l)$,
each unique up to isometry,
as undented hyperbolic polyhedra of side length $l$,
$0 < l \leq \infty$,
where
$\hpuff(\sigma,\infty)$
is an equilateral ideal hyperbolic polyhedron.
After rescaling all lengths by $1/l$,
$\hpuff(\sigma,l)$ converges as $l\goesto0$
to the Euclidean neoplatonic realization $\puff(\sigma)$.
\end{conjecture}

The atlas includes animations showing the deformation from ideal to
Euclidean realizations.

The intuition behind this conjecture is that if we begin with the unique ideal
realization $\hpuff(\sigma,\infty)$ and continue it through a family
$\hpuff(\sigma,l)$, then it is unable to develop a dent, because one cannot
slip a curve whose length remains $< 2\pi$
over the unit sphere.
The issue is therefore whether this homotopy can be continued all the way
down to $l=0$.
If it can, the Euclidean neoplatonic exists.
If in addition there can be no bifurcation along the way, then it is unique.

For every prime $6$-net with $v\leq50$,
we have rigorous existence proofs for both Euclidean and ideal neoplatonics,
but our evidence for the homotopy conjecture is experimental.
The approximate neoplatonics used in the proof of Proposition \ref{euclid}
were obtained by numerical homotopy from ideal neoplatonics.

\section{Realizing ideal neoplatonics} \label{sec:real}

To realize ideal neoplatonics we use the approach of
Bobenko, Pinkall, and Springborn
\cite{bobenko:discrete}
and Springborn
\cite{springborn:ideal}.
Figure \ref{fex} (top) shows a convex ideal hyperbolic polyhedron $M$
viewed through one of its ideal vertices $v_0$.
What one sees is a triangulated Euclidean polygon $P$.
It lives in a visual horosphere about $v_0$,
which can be identified with the sphere at infinity with $v_0$ removed.
The symmetry axes of the ideal triangle faces not incident with $v_0$
appear as the ``symmedians'' of the triangles of $P$, shown in red.
\figsize{.8\linewidth}{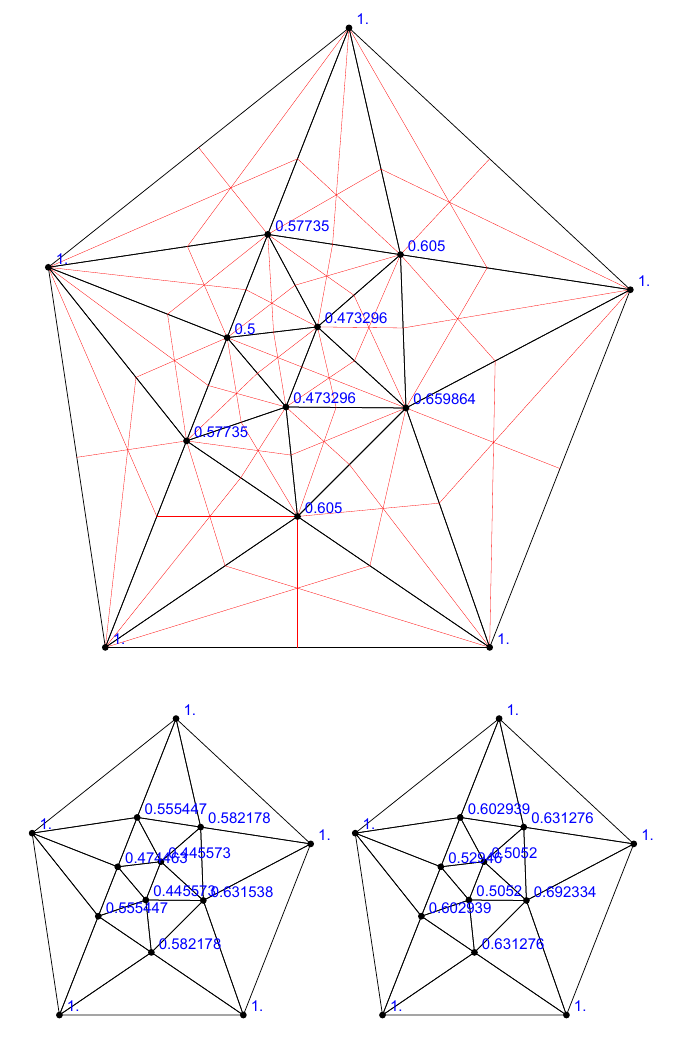}{fex}{Realizing an ideal neoplatonic.}

We can tell that $M$ is convex for two reasons:
\begin{enumerate}
\item
The polygon $P$ is convex;
\item
The triangulation satisfies the Delaunay empty-circle condition:
the circumcircle of each triangle of $P$ contains no vertex of $P$ in its interior.
\end{enumerate}

We can tell that $M$ is equilateral for three reasons:
\begin{enumerate}
\item
The boundary edges of $P$ all have the same length.
This means that there is no shearing along edges incident with $v_0$;
\item
The triangles adjacent to the boundary are isosceles.
This means that there is no shearing along the remaining edge of any
triangle incident with $v_0$;
\item
The red symmedian lines match up where they meet the interior edges of $P$.
This means there is no shearing along edges between triangles not
incident with $v_0$.
\end{enumerate}

There is an algebraic way to express the absence of shearing.
In Figure \ref{fex}
the vertices of $P$ are labeled by a function
$u:V^\prime \to (0,\infty)$.
The edge $ab$ has length
\[
\len(ab)=u(a)u(b).
\]
At the boundary vertices of $P$
(those adjacent to $v_0$ in $M$)
the function $u$ takes the value $1$.
These properties imply
that the quadrilaterals $abcd$ formed by triangles $abc$ and $acd$
adjacent along an interior edge $ac$ are conformal rhombi:
\[
\len(ab)\len(cd) = \len(ad)\len(bc) = u(a)u(b)u(c)u(d)
.
\]
If we think of $P$ as living in the complex plane,
this means that the appropriate cross-ratio of the points $a,b,c,d$ has
modulus $1$.
This is an algebraic way of saying that there is no shearing along $ac$.

As explained by Springborn
\cite{springborn:ideal},
to construct such a realization
for a given triangulation of $S^2$
one seeks an appropriate function $u$
from which to build $P$ by gluing triangles with edge lengths
$\len(ab)=u(a)u(b)$.
We need the following properties:
\begin{enumerate}
\item
For every triangle $abc$, we require the triangle inequalities:
\[
u(a)(u(b)+u(c)) > u(b)u(c)
\]
and its cyclic variants.
\item
The total angle at every interior vertex is $2\pi$,
so that $P$ lies flat.
\item
The total angle at every boundary vertex is $\leq \pi$,
so that the polygon is convex.
\item
The triangulation is Delaunay:
for quadrilaterals $abcd$ as above, writing $U(x)=u(x)^2$, we require
\[
2 U(b)U(d)(U(a)+U(c)) - U(a)U(c)(U(b)+U(d)) \geq 0
.
\]
\end{enumerate}

One possible approach is to
compute $u$ numerically to high precision,
recover from it approximate coordinates for the vertices of $P$,
recognize those coordinates as algebraic numbers,
prove the rhombus equalities algebraically,
and verify the Delaunay inequalities using arithmetic with guaranteed
error bounds.
We ran into trouble recognizing the algebraic
numbers once the number of vertices got bigger than $11$ or so.
Where we succeeded, fields of definition were not special in any
obvious way.

There is a looser approach:
the method of super- and subsolutions,
familiar from the Dirichlet problem and circle packing.
Constructing $P$ resembles circle packing,
except that circle packing assigns the edge $ab$ the length
$u(a)+u(b)$ rather than $u(a)u(b)$
and has no corresponding issue with triangle inequalities.

Let us illustrate this by example.
Figure \ref{fex} (bottom) shows two functions $u_0 \leq u_1$
for an example with $v=14$.
On the boundary $u_0=u_1=1$.
At interior vertices,
$u_1$ has defect $\geq 0$ and $u_0$ has defect $\leq 0$.
(In other words, $u_1$ puts too little angle at interior vertices,
$u_0$ puts too much.)
Starting with $u=u_1$,
we can pick an interior vertex $a$ and adjust the value of $u(a)$
downward to add angle at $a$ and make the defect disappear there.
The defect diffuses to neighboring vertices, whose defects all increase.
We can then go from vertex to vertex,
pushing $u$ down, down, down.
The function $u_0$ acts as a floor.
If $u \geq u_0$, lowering $u(a)$ all the way to $u_0(a)$
would make the defect at $a$ negative.
Thus the defect vanishes before $u(a)$ reaches $u_0(a)$,
and the inequality $u\geq u_0$ is preserved.

What about the triangle inequality?
Could this fail and prevent us from adjusting $u(a)$ to remove the defect?
No, because we can check that as long as $u_0 \leq u \leq u_1$ the
triangle inequalities hold for all triangles of $P$.
So we never run off the rails during the descent.

This produces a solution $u$, bracketed by $u_0$ and $u_1$,
for which $P$ lies flat.
To prove convexity we check that any function bracketed in this way must have
total angle $< \pi$ at boundary vertices,
and satisfy the Delaunay inequality at every interior edge.

This approach works to prove existence of convex ideal neoplatonics for all
prime 6-nets with $4 \leq v \leq 50$.
For $v=4,\ldots,24$,
we chose super- and subsolutions with all defects equal to $\pm1/500$.
For $v=25$ we used $\pm1/1000$;
by $v=50$ we were using $\pm1/4000$.

We think the Perron approach could lead to a general proof.
It reduces the problem to constructing super- and subsolutions
whose entire bracket satisfies the triangle, Delaunay, and boundary inequalities.
The computations show that these bounds must become tight:
at $v=50$ the defects were $\pm1/4000$.

\subsection*{The ideal icosahedron}

We close this section, as Euclid closed Book XIII, with the icosahedron,
only now in its incarnation
as an ideal neoplatonic.
We get this from the functions
$u_0,u_1$ shown at the bottom of Figure \ref{ficos}.
\fig{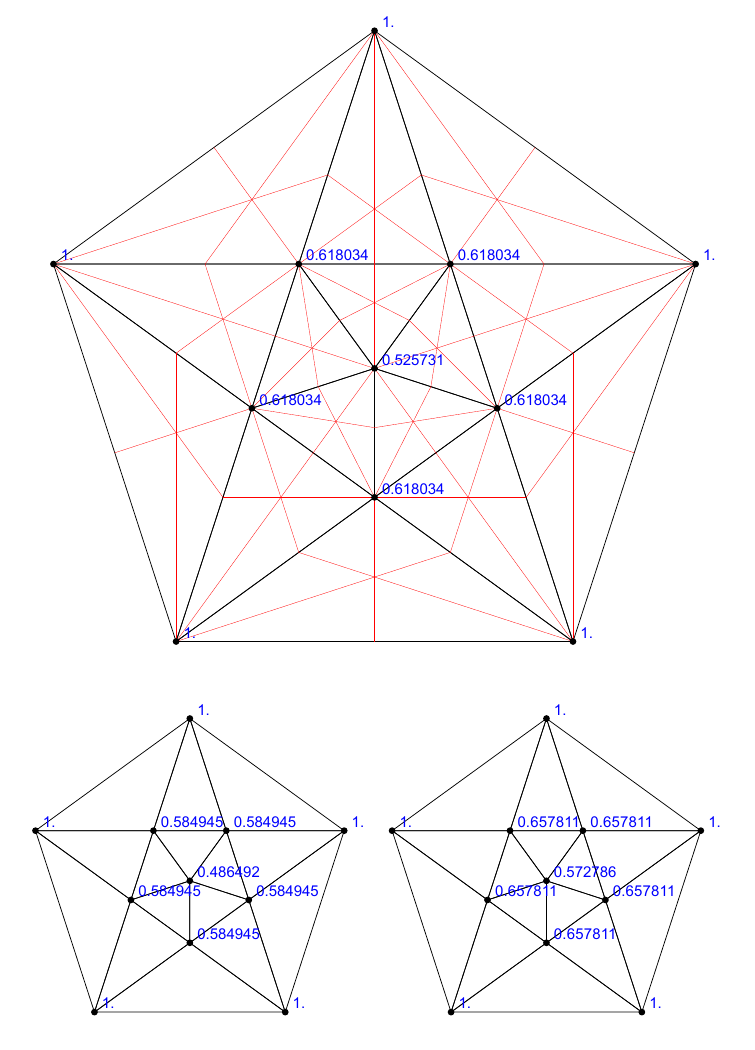}{ficos}{Realizing the ideal icosahedron.}
The exact values for $u$, shown at the top of the figure, are
\[
(\sqrt{5}-1)/2 = .61803\ldots
;\;
\sqrt{(5-\sqrt{5})/10} = .52573\ldots
.\]

A direct verification is also available.
It suffices to check that when we inscribe a pentagram in a regular
pentagon and subdivide the interior pentagon,
every pair of adjacent triangles yields a conformal rhombus.
This is clear by inspection.

Standard methods
(Cf.\ Vinberg
\cite[Sec.\ 5.2]{vinberg:volume})
give the volume of the ideal icosahedron as
\[
13.529628198483504531786905296849240634601965635662\ldots
.
\]
We do not know a standard reference in which this value is tabulated.

\section{Further conjectures}

\subsection*{Neoconvex bodies}

If neoplatonics are unique, it is natural to ask whether
uniqueness holds more generally.

Call a hyperbolic or Euclidean polyhedron
\emph{neoconvex} if it is
\emph{intrinsically convex} (all angle defects are $\geq 0$) and undented.
This is a discrete analogue of having nonnegative Gaussian curvature
and nonnegative mean curvature,
which for a smooth surface would be equivalent to convexity.

\begin{hypothesis}[Neoconvex rigidity]
Two neoconvex polyhedra with the same combinatorics and isometric faces
are ambiently isometric.
\end{hypothesis}

This would imply uniqueness of neoplatonics.
If true, neoconvex rigidity would
strengthen Cauchy rigidity
\cite{cauchy:rigid}
in a direction orthogonal
to Alexandrov-type uniqueness
\cite[Lemma 1, p.\ 62]{alex:convex}.
Alexandrov allows faces to bend while insisting on global convexity.
Here the faces remain rigid, but nonconvexity is allowed,
provided that the vertices remain undented.
Even an ordinary convex polyhedron would then be unique
among all undented polyhedra with the same metric net,
not merely among convex ones.

A stronger version is:

\begin{hypothesis}[Dent location]
An embedded polyhedron with no vertices of negative angle defect is unique up to isometry
among embedded realizations of its metric net with the same set of dented vertices.
\end{hypothesis}

Embeddedness is essential:
denting a vertex of the regular octahedron gives the flexible
``cootie catcher,'' which can fold onto a single triangle.

Under this hypothesis,
a polyhedron with no vertex of negative angle defect
could at most have $2^v$ embedded realizations up to isometry.
Thus a continuous flex through embedded polyhedra with no vertices of negative angle defect
would disprove the hypothesis.

\section{Acknowledgments and computations}

These conjectures emerged from playing with models built from
Polydrons\texttrademark\ and Geometiles\texttrademark,
in the context of our work with Zili Wang on the conjecture of
Sleator, Tarjan, and Thurston \cite{stt:jams}.
See, for example, \cite[Fig.\ 15]{dew:filling}.
We are deeply indebted to Zili for her support and encouragement
throughout.

To construct and prove the existence of ideal realizations,
we used the approach of Springborn
\cite{springborn:ideal},
which builds on work of Luo \cite{luo:yamabe}
and Bobenko, Pinkall, and Springborn
\cite{bobenko:discrete}.

Further inspiration came from
Agol
\cite[`frog eggs', p.\ 5]{agol:haken};
Sabitov
\cite{sabitov:volume};
Connelly
\cite{connelly:alexandrov};
Audet
\cite{audet:cayleymenger};
and extensive discussions with Yana Mohanty and Richard Schwartz.

Computations were carried out primarily by Claude Code,
with supervision and advice from ChatGPT.
Both contributed to the development of ideas, and the exposition benefited
from ChatGPT.
The proof programs and certificate formats are distributed on GitHub in the repositories
\href{https://github.com/peterdoyle1717/bendprover}{\nolinkurl{peterdoyle1717/bendprover}} and
\href{https://github.com/peterdoyle1717/ideal}{\nolinkurl{peterdoyle1717/ideal}}.

We encoded nets using the edgebreaker algorithm of
Rossignac, Safonova, and Szymczak
\cite{edgebreaker, clers}.
We assigned each net a canonical CLERS name
by changing the edgebreaker code letters from CLERS to CBEAD
and taking the lexicographically first representative.

We generated fullerene duals using the buckygen software of
Brinkmann, Goedgebeur, and McKay
\cite{buckygen},
and promoted them to a catalog of prime $6$-nets
by a variant of the method of Goedgebeur
\cite{ghent},
using the canonical CLERS name to identify isomorphic nets.

We are grateful to all these researchers for sharing their
ideas, algorithms, and code.

\clearpage
\begin{appendices}
\section{Counts of prime 6-nets}\label{primecount}
The table gives the numbers $p_v$ of prime 6-nets with $v$
vertices, $1\le v\le 100$. They were computed using a variant of the
method of Goedgebeur \cite{ghent}, and can be deduced from the tables presented there.

\begin{center}
  \small
  \setlength{\tabcolsep}{3pt}
  \begin{tabular}{r@{\hspace{0.45em}}l|r@{\hspace{0.45em}}l|r@{\hspace{0.45em}}l|r@{\hspace{0.45em}}l}
  \toprule
  $v$ & $p_v$ & $v$ & $p_v$ & $v$ & $p_v$ & $v$ & $p_v$ \\
  \midrule
   1 & 0     & 26 & 6\,539        & 51 & 1\,651\,488   & 76 & 47\,782\,118  \\
   2 & 0     & 27 & 8\,750        & 52 & 1\,945\,387   & 77 & 53\,374\,432  \\
   3 & 0     & 28 & 11\,854       & 53 & 2\,279\,167   & 78 & 59\,595\,319  \\
   4 & 1     & 29 & 15\,629       & 54 & 2\,668\,297   & 79 & 66\,400\,484  \\
   5 & 0     & 30 & 20\,755       & 55 & 3\,109\,939   & 80 & 73\,920\,746  \\
   6 & 1     & 31 & 26\,924       & 56 & 3\,619\,603   & 81 & 82\,155\,286  \\
   7 & 1     & 32 & 35\,070       & 57 & 4\,196\,607   & 82 & 91\,253\,047  \\
   8 & 2     & 33 & 44\,866       & 58 & 4\,862\,665   & 83 & 101\,140\,713 \\
   9 & 3     & 34 & 57\,649       & 59 & 5\,608\,787   & 84 & 112\,064\,305 \\
  10 & 7     & 35 & 72\,749       & 60 & 6\,464\,974   & 85 & 123\,946\,248 \\
  11 & 10    & 36 & 91\,995       & 61 & 7\,428\,133   & 86 & 136\,995\,975 \\
  12 & 22    & 37 & 114\,966      & 62 & 8\,523\,556   & 87 & 151\,178\,908 \\
  13 & 32    & 38 & 143\,725      & 63 & 9\,748\,741   & 88 & 166\,745\,838 \\
  14 & 59    & 39 & 177\,439      & 64 & 11\,142\,291  & 89 & 183\,587\,072 \\
  15 & 92    & 40 & 219\,402      & 65 & 12\,692\,741  & 90 & 202\,085\,816 \\
  16 & 153   & 41 & 268\,496      & 66 & 14\,449\,523  & 91 & 222\,059\,220 \\
  17 & 230   & 42 & 328\,571      & 67 & 16\,402\,234  & 92 & 243\,882\,591 \\
  18 & 373   & 43 & 398\,578      & 68 & 18\,597\,312  & 93 & 267\,490\,334 \\
  19 & 536   & 44 & 483\,170      & 69 & 21\,036\,300  & 94 & 293\,261\,943 \\
  20 & 820   & 45 & 581\,293      & 70 & 23\,778\,320  & 95 & 320\,966\,986 \\
  21 & 1\,174 & 46 & 699\,570     & 71 & 26\,798\,812  & 96 & 351\,250\,688 \\
  22 & 1\,735 & 47 & 835\,175     & 72 & 30\,186\,590  & 97 & 383\,795\,613 \\
  23 & 2\,418 & 48 & 996\,468     & 73 & 33\,926\,502  & 98 & 419\,202\,743 \\
  24 & 3\,442 & 49 & 1\,182\,607  & 74 & 38\,091\,396  & 99 & 457\,264\,253 \\
  25 & 4\,711 & 50 & 1\,401\,622  & 75 & 42\,676\,406  & 100 & 498\,626\,241 \\
  \bottomrule
  \end{tabular}
\end{center}
\end{appendices}

\clearpage
\bibliography{neo}
\bibliographystyle{amsplain}
\end{document}